\documentclass{amsart}

\usepackage{amsmath}
\usepackage{amssymb}
\usepackage{eucal}
\usepackage{amsthm}
\usepackage{verbatim}
\usepackage[hyphens]{url}
\usepackage{natbib}

{\theoremstyle{definition}

}

{\theoremstyle{remark}

}

\DeclareMathOperator{\Prob}  {\mathsf{P}}
\DeclareMathOperator{\Expec} {\mathsf{E}}

\newcommand{\bsone}{\boldsymbol{1}}

\newcommand{\bsl}{\boldsymbol{l}}

\newcommand{\mcJ}{\mathcal{J}}
\newcommand{\smcJ}{ \text{\raisebox{0.15ex}{\small$\mathcal{J}$}} }

\newcommand{\mcL}{\mathcal{L}}

\newcommand{\mcV}{\mathcal{V}}

\newcommand{\mbR}{\mathbb{R}}

\newcommand{\stfrac}[2]{ \leavevmode\ensuremath{
    \kern.1em\raise.5ex\hbox{\footnotesize #1}
    \kern-.1em/
    \kern-.15em\lower.25ex\hbox{\footnotesize #2}} }
\newcommand{\smfrac}[2]{ \leavevmode\ensuremath{
    \kern.1em\raise.5ex\hbox{\footnotesize $#1$}
    \kern-.1em/
    \kern-.15em\lower.25ex\hbox{\footnotesize $#2$}} }
\newcommand{\Stfrac}[2]{ \leavevmode\ensuremath{
    \kern.1em\raise.5ex\hbox{#1}
    \kern-.1em/
    \kern-.15em\lower.25ex\hbox{#2}} }
\newcommand{\Smfrac}[2]{ \leavevmode\ensuremath{
    \kern.1em\raise.5ex\hbox{$#1$}
    \kern-.25em/
    \kern-.15em\lower.5ex\hbox{$#2$}} }

\begin{document}

\title{A New Efficient Algorithm
       for Construction of LLS Models}

\author[M.Kovtun]{Mikhail~Kovtun}
\address{Duke University, CDS\\
         2117 Campus Dr., Durham, NC, 27708}
\email{mkovtun@cds.duke.edu}

\author[I.Akushevich]{Igor~Akushevich}
\address{Duke University, CDS\\
         2117 Campus Dr., Durham, NC, 27708}
\email{aku@cds.duke.edu}

\author[K.G.Manton]{Kenneth~G.~Manton}
\address{Duke University, CDS\\
         2117 Campus Dr., Durham, NC, 27708}
\email{kgm@cds.duke.edu}

\author[H.D.Tolley]{H.~Dennis~Tolley}
\address{Brigham Young University, Department of Statistics\\
          206 TMCB, Provo, UT}
\email{tolley@byu.edu}

\thanks{This research was supported by grants
        from National Institute of Aging.}

\thanks{This work was presented on Annual Meeting
        of Population Association of America,
        Philadelphia, PA, April 2005.}

\subjclass[2000]{Primary 62H12; Secondary 62J99}

\begin{abstract} 
We present a new efficient algortithm for construction
of linear latent structure (LLS) models.
This algorithm reduces a problem of estimation of model parameters
to a sequence of problems of linear algebra,
which assures a low computational complexity
and ability to handle on desktop computers data that involve
up to thousands of variables.
\end{abstract} 

\maketitle 

\bigskip

The class of linear latent structure (LLS) models
belongs to a family of latent structure models,
which, in turn, is a subfamily of a family
of mixed distribution models. Such models naturally occur
when a population of interest is supposed to be heterogeneous.

The most widely used methods for estimation latent structure
models are based on maximization of the likelihood function.
These are well established methods possessing many good properties.
Nevertheless, they have limitations, which may restrict
or even prevent their usage.
First, the number of parameters to be optimized
is proportional to the number of variables (measurements),
which in practice limits the number of variables used in
the analysis to several dozens.
Second, the likelihood function in the case of latent structure
analysis is often multimodal, which requires usage of additional
techniques to ensure that the absolute maximum is found.

Our algorithm is based on methods of linear algebra,
which eliminates the problem
of multimodality and allows us to analyze
up to thousands of variables.
The time spent by the algorithm is proportional to
the cube of the number of variables.

Historically, the predecessor of LLS analysis was
grade of membership (GoM) analysis, which was introduced in 
\cite{Woodbury:1974}; see also \cite{Manton:1994} for detailed
exposition and additional references.
Our work on LLS analysis originated from attempts to find
conditions for consistency of GoM estimators.
The development eventually lead to a new class of models,
which differ from GoM models in a way how the model is
formulated, methods of model estimation, meaning of estimators
and their interpretation.

\section{Basic notions}

LLS analysis considers $J$ discrete measurements,
represented by a vector of random variables $X = (X_1, \dots, X_J)$,
with the set of outcomes of $j^\text{th}$ measurement
(i.e. the set of possible values of random variable $X_j$)
being $\{1, \dots, L_j\}$.
We consider a distribution law of this random vector
as a mixture of {\em independent} distribution laws,
i.e. distribution laws satisfying

\begin{equation}
\label{eq:Independency}
\Prob( X_1=l_1 \wedge \dots \wedge X_J=l_J )
= \prod_j \Prob( X_j=l_j )
\end{equation}

Representation of the observed distribution law as a mixture
of independent distribution laws is standard for latent structure
analysis (and it is its defining characteristic).

Due to (\ref{eq:Independency}), description of independent
distribution law requires only knowing $\Prob(X_j = l_j)$.
Thus, an independent distribution law may be described by
$|L|=L_1+\dots+L_J$-dimensional vector
$\beta=(\beta_{jl})_{jl}$,
where $j$ ranges from $1$ to $J$, and for every $j$,
$l$ ranges from $1$ to $L_j$; $\beta_{jl} = \Prob(X_j = l)$.
Let $\mu_\beta$ be a mixing measure producing the observed
distribution.

The main LLS assumption is that for some integer $K$,
$\mu_\beta$ is supported by a $K$-dimensional
linear subspace $Q$ of $\mbR^{|L|}$.
Later, we refer to this $K$ as to the
{\em dimensionality of LLS problem}.

This assumption is essentially equivalent to the assumption that
there exists a $K$-di\-men\-si\-o\-nal random vector $G$
such that for every $j$ a regression of $Y_j$ on $G$ is linear.
Here $Y_j$ is an $L_j$-dimensional random vector,
$Y_j = \bsone_l$ if $X_j = l$
(where $\bsone_l$ denotes a vector which has $l^\text{th}$ component
equal to $1$, and all other components equal to $0$.)
Namely, let $\Lambda=\{\lambda^1,\dots,\lambda^K\}$,
$\lambda^k=(\lambda^k_{jl})_{jl}$, be any basis of $Q$,
and for $\beta \in Q$, let $g=(g_k)_{k=1,\dots,K}$ be
its coordinates in basis $\Lambda$. Then the random vector $G$
is the random vector $\beta$ (distributed according $\mu_\beta$)
written in coordinates $g$, and matrices
$\Lambda_j = (\lambda^k_{jl})_{kl}$ are linear regression matrices.

The linear regression assumption is crucial for understanding
the meaning of the LLS model and gives guidelines for its
applicability.
It essentially means that the measurements are not chosen
arbitrarily but rather to reflect in some degree a hidden
property, or a hidden state, represented by the random vector $G$.
LLS analysis is about how to discover this hidden state
and describe it as precisely as possible.

Let $\mu_g$ be a measure $\mu_\beta$ written in coordinates $g$.

Let $\ell = (\ell_1,\dots,\ell_J)$ be an integer vector with
$0 \le \ell_j \le L_j$.
Such a vector represents the outcome of $J$ measurements,
and $\ell_j=0$ means that we do not take into account
the outcome of the $j^\text{th}$ measurement.
Thus, a value of $\ell_j = 0$ in a vector $\ell$ means that the
vector is a marginal vector across all values of the $j^\text{th}$
measurement.
Let $\mcL^0$ be a set of all such vectors,
and for every $\smcJ \subseteq \{1,\dots,J\}$ let $\mcL^{[\mcJ]}$
be a set of vectors having $0$'s exactly on places from $\smcJ$.
Let $v = (v_1,\dots,v_K)$ be an integer vector with $v_k \ge 0$,
and for every integer $J' \ge 0$ let $\mcV[J']$ be a set of such vectors
satisfying the additional condition $\sum_k v_k = J'$.

In this language, the values of interest are unconditional moments
of the distribution $\mu_\beta$

\begin{equation}
M_\ell(\mu_\beta) =
    \int \prod_{j \,:\, \ell_j \neq 0} \beta_{j \ell_j}
    \, \mu_\beta(d\beta)
\end{equation}

\noindent
and conditional moments of distribution $\mu_g$,

\begin{equation}
\label{eq:CondMoments}
\Expec(G^v \mid X=\ell) =
    \int \prod_k g^{v_k}_k
    \frac{\prod_{j \,:\, \ell_j \neq 0} \sum_k g_k \lambda^k_{jl}}
         {M_\ell(\mu_\beta)}
    \, \mu_g(dg)
\end{equation}

The unconditional moments $M_\ell(\mu_\beta)$ are the probabilities
of obtaining the response pattern $\ell$
(under assumptions of the model.)
Thus, frequencies of response patterns $\ell$ in a sample,
denoted $f_\ell$,
are consistent and efficient estimators for unconditional moments
$M_\ell(\mu_\beta)$.

The conditional moments $\Expec(G^v \mid X=\ell)$ express
our knowledge of the state of the individual (represented by
random vector $G$) based on the outcomes of the measurements.
These values are not directly estimable from the observations.
The goal of LLS analysis is to obtain estimates for these
conditional moments.

The most important relation connecting unconditional moments,
conditional moments and the basis $\Lambda$ (in which conditional
moments are calculated) is:

\begin{equation}
\label{eq:MainEq}
\sum_k \lambda^k_{jl} \cdot
\left(
    M_\ell(\mu_\beta) \cdot
    \Expec(G^{v+\bsone_k} \mid X=\ell)
\right) =
M_{\ell+\bsl_j}(\mu_\beta) \cdot \Expec(G^v \mid X=\ell+\bsl_j)
\end{equation}

\section{The main system of equations}

We have shown in \cite{Kovtun:2004} that the LLS model defined above
is fully described by a system
of equations (with respect to variables $\alpha^k_{jl}$ and $h^v_\ell$)

\begin{equation}
\label{eq:MainEqSys}
\begin{cases}
\sum_k \alpha^k_{jl} h^{v+\bsone_k}_\ell = h^v_{\ell+\bsl_j}, \quad
    & J' \in [0..J-1],~~
      v \in \mcV[J'],\\
    & \smcJ \subseteq [1..J] \,:\, |\smcJ|>J',~~
      \ell \in \mcL^{[\mcJ]}_{\phantom{L}},\\
    & j \in \smcJ,~~
      l \in [1..L_j]
\\[3pt]
h^{(0,\dots,0)}_\ell = M_\ell,
    & \ell \in \mcL^0
\\[3pt]
\sum_{v \in \mcV[J']}
        \frac{(\sum_k v_k)!}{\prod_k v_k!} h^v_{(0,\dots,0)} = 1,
    & J' \in [0..J]
\end{cases}
\end{equation}

\noindent
In this system, the first group of equations corresponds to the
main relation between moments (\ref{eq:MainEq}), and the last
two equations are normalization conditions.

We have proven the following properties of the main system:

\begin{enumerate}
\item
    Any basis $\Lambda$ of $Q$ together with conditional moments
    $\Expec(G^v \mid X=\ell)$ calculated in this basis give a
    solution of (\ref{eq:MainEqSys})
    ($\lambda^k_{jl}$ should be substituted for $\alpha^k_{jl}$, and
    $M_\ell(\mu_\beta) \cdot \Expec(G^v \mid X=\ell)$
    should be substituted for $h^v_\ell$.)
\item
    Under mild conditions, {\em every} solution of (\ref{eq:MainEqSys})
    gives a basis of $Q$ and conditional moments calculated in this
    basis.
\end{enumerate}

As the main system of equations fully describes the model,
the important property of the LLS analysis follows:
the mixing distribution is not fully identifiable.
Only a finite number of moments may be found by solving the system,
and any mixing distribution that have these moments would satisfy
the system.
The fact of nonidentifiability also follows from the general
theorem about identifiability of mixtures,
because the family of distributions contained in $Q$
is not linearly independent.

The attractive feature of the LLS analysis is that
it can discover a number of useful invariants of the mixture.
The supporting plane of the mixing distribution is defined
uniquely, and low-order moments are identifiable as well.
This information is sufficient to make practically substantial
conclusions about the population under consideration.

The main system of equations provides a means for consistent
estimation of model parameters. The solution of this
system continuously depends on unconditional moments $M_\ell$;
thus, substitution of frequencies $f_\ell$ for moments $M_\ell$
gives a system, which solutions converge to the true values
of parameters when frequencies converge to the true moments.

One good property of the main system of equations
is that it is linear with respect to variable $h^v_\ell$.
Thus, if the supporting plane of distribution is known,
the conditional moments (\ref{eq:CondMoments})
may be obtained by solving a linear
system of equations. It happens that the supporting plane
may be found independently by analysis of the moment
matrix, which we describe in the next subsection.

\section{The moment matrix}

Let us write a vector of moments $(M_{\bsl_j})_{jl}$ together with
incomplete vectors of moments
$(M_{\bsl'_{j'}+\bsl_j})_{jl \,:\, j \neq j'}$,
etc.,
as columns of
a matrix, with places for which we do not have moments filled
by question marks. We refer to this incomplete matrix as the
{\em moment matrix}.
The moment matrix contains a column for every
$\ell \in \mcL^0$.
Figure \ref{fig:1} gives an example of a portion of a moment matrix for
the case $J=3$, $L_1=L_2=L_3=2$.
Columns in this matrix correspond to
$\ell=(000)$, $(100)$, $(200)$, $(010)$, $(020)$, $(001)$, $(002)$,
$(110)$;
other columns are not shown.

\begin{figure}[ht]
\begin{equation*}
\begin{pmatrix}
M_{(100)}   &
    ?           &  ?            &
    M_{(110)}   &  M_{(120)}    &
    M_{(101)}   &  M_{(102)}    &
    ?           &  \cdots \phantom{\vdots}       \\
M_{(200)}   &
    ?           &  ?            &
    M_{(210)}   &  M_{(220)}    &
    M_{(201)}   &  M_{(202)}    &
    ?           &  \cdots \phantom{\vdots}       \\
M_{(010)}   &
    M_{(110)}   &  M_{(210)}    &
    ?           &  ?            &
    M_{(011)}   &  M_{(012)}    &
    ?           &  \cdots \phantom{\vdots}       \\
M_{(020)}   &
    M_{(120)}   &  M_{(220)}    &
    ?           &  ?            &
    M_{(021)}   &  M_{(022)}    &
    ?           &  \cdots \phantom{\vdots}       \\
M_{(001)}   &
    M_{(101)}   &  M_{(201)}    &
    M_{(011)}   &  M_{(021)}    &
    ?           &  ?            &
    M_{(111)}   &  \cdots \phantom{\vdots}       \\
M_{(002)}   &
    M_{(102)}   &  M_{(202)}    &
    M_{(012)}   &  M_{(022)}    &
    ?           &  ?            &
    M_{(112)}   &  \cdots \phantom{\vdots}      \\[-6pt]
&
\end{pmatrix}
\end{equation*}
\caption{\label{fig:1} Example of moment matrix}
\end{figure}

Note that certain moments (which are replaced by question marks
in the moment matrix) are not observable.
The reason for this is that we are not able
to perform a measurement on an individual multiple times independently,
and since individuals are heterogeneous (have different probabilities
of outcomes of measurements), we do not have multiple realizations
of independent identically distributed random variables.

For a moment matrix $M$ let its completion $\bar{M}$ be a matrix
obtained from $M$ by replacing question marks by arbitrary numbers.
We have shown that the moment matrix always has a completion
in which all columns belong to the supporting plane $Q$.
Thus, if the moment matrix has sufficient rank
(which is the case in non-degenerate situations,)
a basis of $Q$ may be obtained from this matrix.
As we have a consistent estimator of the moment matrix in form
of a frequency matrix, the supporting plane may be
consistently estimated.

This property of the moment matrix suggests an efficient
algorithm to obtain LLS estimates.
First, a basis of the supporting plane can be obtained from
the moment matrix (a way to do this is described in the next section),
and second, conditional moments can be found by solving
a linear system of equations.

\section{Algorithm}
\label{ssec:Algorithm}

As it is suggested by a structure of the main system of equations
(\ref{eq:MainEqSys}) and by properties of the moment matrix,
the algorithm is naturally decomposed into two parts.
On the first step, a basis of the supporting plane
should be constructed; the input for this step is the frequency matrix.
On the second step, a system of linear equations should be solved
to obtain estimates for conditional expectations.

\paragraph{Step 1: Finding the supporting plane.}
As for the model distribution all columns of the moment matrix
belong to the supporting plane, and as the frequency matrix
is an approximation of the moment matrix,
the natural way to search for the supporting plane is
to search for a plane that minimizes the sum of distances
from it to the columns of the frequency matrix.
In our case, however, this way is complicated by
at least three obstacles:
(a) a sought basis $\Lambda$ should {\em exactly} satisfy
conditions $\sum_l \lambda^k_{jl} = 1$ for every $k$ and $j$;
(b) the statistical inaccuracy of approximation of moments $M_\ell$ by
frequencies $f_\ell$ varies considerably over elements of
frequency matrix;
(c) the moment matrix (and, correspondingly, the frequency matrix)
is incomplete.

The suggested algorithm for estimating the supporting plain
consists of the following steps.

{\renewcommand{\theenumi}{\roman{enumi}}
 \renewcommand{\labelenumi}{(\theenumi)}
\begin{enumerate}
\item
    The computational rank of the frequency matrix is estimated.
    For this, we take the biggest minor of the frequency matrix
    that does not contain question marks. (For the example given
    in Figure \ref{fig:1}, it is the left bottom minor of size
    $3 \times 3$.) Then we calculate the singular value decomposition
    (SVD) and take $K_0$ (the first approximation of dimensionality
    of the LLS problem) equal to the number of singular values
    that are greater than standard deviation of the norm of columns
    involved in the minor.
    (The final value for dimensionality of LLS problem will be
    chosen on the step (v).)
    As one of requirements for applicability of LLS model is
    $K \ll |L|$, nothing can be done further if all (or too many)
    singular values are greater than the standard deviations.
\item
    We construct a completion of the frequency matrix by means
    of the following procedure. For every column $c$ of the
    frequency matrix and row $jl$ of a question mark in $c$,
    we select $K_0$ columns $c^{(1)},\dots,c^{(K_0)}$ satisfying:
    (a) all columns $c^{(i)}$ contain a value (not a question mark)
    in row $jl$; (b) there exist $p \ge K_0$ rows such that all
    columns $c, c^{(1)},\dots,c^{(K_0)}$ contain values in these
    rows. Let $c[p]$ be a subcolumn containing only selected rows.
    Then we solve a linear system
    $\alpha_1 c^{(1)}[p] + \dots + \alpha_1 c^{(K_0)}[p] = c[p]$
    and replace a question mark at the position $c_{jl}$ by
    $\alpha_1 c^{(1)}_{jl} + \dots + \alpha_1 c^{(K_0)}_{jl}$.
    The system to be solved is overdetermined; we solve it
    by minimization of residuals using SVD.
    When $K_0$ is sufficiently smaller than $|L|$, the required
    selection of the columns is possible for every column $c$
    which contains at least $K_0$ values; columns containing
    less than $K_0$ values are discarded from further consideration.
    According to \cite{Kovtun:2004}, the moment matrix always
    has a completion of rank equal to the dimensionality of LLS
    problem; thus, this method should give good results.
\item
    Columns are normalized, so that the condition
    $\sum_l c'_{jl}=1$ holds for every $j$.
    This is always possible, as for every column $c$ we have
    $\sum_l c_{jl} = s$, where $s$ does not depend on $j$.
    Thus, we take $c' = \frac{1}{s} c$.
\item
    Next, we remove the restriction $\sum_l c_{jl} = 1$ by
    reducing number of rows by $J$ ($1$ for every group of indices
    $j1,\dots,jL_j$). For this, we use a linear map
    from $\mbR^{|L|}$ to $\mbR^{|L|-J}$ given by a block-diagonal
    matrix $A$ with $J$ blocks
    \begin{equation}
    \label{eq:RemovingDependencies}
        A_j = \begin{pmatrix}
                - \frac{\sqrt{L_j}-1}{L_j-1}   & 1 & 0 & \dots & 0 \\
                  \hdotsfor{5}                                     \\
                - \frac{\sqrt{L_j}-1}{L_j-1}   & 0 & 0 & \dots & 1
              \end{pmatrix}
    \end{equation}
    of size $L_j \times (L_j-1)$.
    This map is an isometry of the subspace of $\mbR^{|L|}$
    defined by equations $\sum_l c_{jl}=1$ to $\mbR^{|L|-J}$
    (every block $A_j$ defines a rotation of a unit simplex in
    $L_j$-dimensional space around hypersurface opposite to the
    first vertex; the angle of this rotation is such that
    the first vertex moves to the point with the first coordinate
    equals 0).
\item
    Now we have $n$ points $y^1, \dots, y^n$
    (images of columns of frequency matrix)
    in $m=|L|-J$-dimensional space. The problem is to find an affine
    plane that minimally deviates from these points. First, we find
    the center of gravity of this system
    \begin{equation}
    \label{eq:CenterOfGravity}
        y^0 = \frac{1}{n} \sum_i y^i
    \end{equation}
    and then consider a new set of points $x^i = y^i-y^0$.
    We need to find a linear subspace in $\mbR^n$ that minimally
    deviates from this set of points. The solution of this problem
    is well-known (see, for example, chapter 43 of \cite{Kendall:1977}):
    one has to consider an $m \times m$  matrix
    $X = (X_{rs})_{rs}$ with components
    $X_{rs} = \sum_i x^i_r \cdot x^i_s$;
    this matrix is symmetrical and positively defined, and thus
    it possess an orthonormal basis of eigenvectors.
    Let $\gamma_1 \ge \gamma_2 \ge \dots \ge \gamma_m \ge 0$ be
    eigenvalues of matrix $X$, and let $z^1,\dots,z^m$ be
    corresponding them eigenvectors.
    The plane of dimensionality $p$ that minimizes the sum of
    squared distances from point $x^1, \dots, x^n$ is spanned by
    $z^1,\dots,z^p$, and the sum of squared distances is
    $\text{tr}X - \sum_{k=1}^p \gamma_k$,
    This gives us a criterion for the selection of the dimensionality
    $K$ of the LLS problem: one has to take $K$ to be the smallest
    integer such that eigenvalues $\gamma_{K}, \dots, \gamma_m$
    are smaller that inaccuracy in input data.
    Vectors $y^0, y^0+z^1, \dots, y^0+z^{K-1}$
    give us an affine basis of the sought affine plane.
\item
    Lastly, we apply inverses of transformation
    (\ref{eq:RemovingDependencies})
    to $y^0, y^0+z^1, \dots, y^0+z^{K-1}$ to obtain the sought basis
    $\lambda^1, \dots, \lambda^K$ of the subspace $Q$.
\end{enumerate}}

The above algorithm solves the problems (a)--(c) listed in the
beginning of the subsection,
and it possesses two important properties that
are crucial to its usefulness:
(a) if the input of the algorithm are true moments of a
distribution generated by $K$-dimensional LLS model, the
output of the algorithm is the true supporting plane;
(b) there exists an open neighbourhood of the true moment
matrix in which the output of the algorithm continuously
depend on its input.

The preliminary experiments with the prototype of the algorithm
performed by the applicants demonstrated that it restores
the supporting plane with a good degree of precision.

\paragraph{Step 2: Calculation of conditional expectations.}
When a basis of the supporting plane is found, the conditional
expectations can be found from the main system of equations
(\ref{eq:MainEqSys}), which is a linear system after substituting
the basis.
This is a sparse overdetermined system; methods for solving
such systems are well-elaborated---see, for example,
\cite{Forsythe:1977, Kahaner:1988}.


\hphantom{
\cite{Everitt:1981}
\cite{Titterington:1985}
}
\hphantom{
\cite{Lindsay:1995}
\cite{Lazarsfeld:1968}
}
\hphantom{
\cite{Bartholomew:1999}
\cite{Bishop:1975}
}
\hphantom{
\cite{Goodman:1978}
\cite{Heinen:1996}
}
\hphantom{
\cite{Langeheine:1988}
\cite{Manton:1994}
}
\hphantom{
\cite{Kovtun:2004b}
\cite{Woodbury:1974}
\cite{Woodbury:1978}
}
\vspace{-1.0in}


\bibliographystyle{apalike}
\bibliography{LLS-PAA2005}

\end{document}